\documentclass{article}

\usepackage{arxiv}

\usepackage{amsmath, cite, amssymb, amsthm}

\usepackage[utf8]{inputenc} 
\usepackage[T1]{fontenc}    

\usepackage{hyperref}       
\usepackage{url}            
\usepackage{booktabs}       
\usepackage{amsfonts}       
\usepackage{nicefrac}       
\usepackage{microtype}      
\usepackage{lipsum}
\usepackage{graphicx}
\usepackage{xcolor}

\title{ON THE CONTROLLABILITY OF AN ORBITING SATELLITE MODEL WITH ELECTROMAGNETIC-ONLY ACTUATION}

\author{
 Ye.~Yevgenieva \\
  Institute of Applied Mathematics and Mechanics, Sloviansk, Ukraine, \\
  Max Planck Institute for Dynamics of Complex Technical Systems, Magdeburg, Germany \\
  \texttt{yevgeniia.yevgenieva@gmail.com} \\
   \And
 A.~Zuyev \\
  Institute of Applied Mathematics and Mechanics, Sloviansk, Ukraine, \\
  Max Planck Institute for Dynamics of Complex Technical Systems, Magdeburg, Germany \\
  \texttt{alexander.zuyev@gmail.com} \\
  \And
 Ju.~Kalosha \\
  Institute of Applied Mathematics and Mechanics, Sloviansk, Ukraine, \\
  \texttt{julykucher@gmail.com} \\
}

\begin{document}
\maketitle

\centerline{{\sc This paper is dedicated to the $100^{\text{th}}$ anniversary of Professor Pavel~V.~Kharlamov.}}\vspace{4ex}

  \centerline
  {\large \bfseries \scshape Abstract}
This paper presents sufficient conditions for small-time local controllability of a control-affine system that describes the rotational motion of a satellite in a circular orbit.
The satellite is modeled as a rigid body subject to electromagnetic actuation.
We focus on the underactuated scenario where the control torque is generated solely by magnetorquers.
The main contributions of this work include proving small-time local controllability around the relative equilibrium under some natural assumptions on the mass distribution of the rigid body.
This result is based on the Lie algebra rank condition and Sussmann's controllability condition.
Furthermore, it is shown that the linearized system is not controllable in a neighborhood of the considered equilibrium.
\vspace{1mm}


\section{Introduction}
One of the fundamental problems in control theory is controllability: determining whether it is possible to drive the system from one state to another using a specified class of controls.
Despite significant advances in this area, the study of controllability for general nonlinear systems remains a challenging problem~\cite{Ag1999,JZC2023,GLPM2024}.
Powerful methods for the controllability analysis of essentially nonlinear systems are extensively covered in~\cite{Cor2007}.
A classical approach relies on the Lie algebraic methods~\cite{AgSa2013}.
For driftless control-affine systems, the controllabity conditions are comple\-te\-ly characterized by the Chow--Rashevskii theorem~\cite{Chow1940,Ra1938}.
For wide classes of nonlinear systems with drift, some important sufficient controllability conditions have been developed in~\cite{BLL2000,Ka1990,He1978,Sus1983,Ste1986,Sus1987,SuJu1972a,SuJu1972b,He1974,He1982}.

This paper is devoted to the mathematical model of a rotating satellite with electro\-magnetic actuation.
To our knowledge, there are no explicit controllability conditions available for this model as of now.
It is mentioned in~\cite{Wi1996} that for a satellite controlled by a set of magnetorquers, magnetic torquing lacks controllability in the direction of the local geomagnetic field.
The controllability problem for a linear model of a spacecraft actuated by magnetic torques is addressed in~\cite{Ya2016}.

In~\cite{Cr1984,AgSa2013}, the Lie algebraic methods have been applied to describe reachability sets and controllability conditions for rotating rigid bodies.
It is known that these kinds of systems also present significant challenges to the stabilization problem.
Notably, it was shown in~\cite{By2008} that a mathematical model of an underactuated rigid spacecraft does not satisfy Brockett's necessary stabilizability condition.

The present work is motivated by the study~\cite{MWZ2022} on a satellite model with electro\-magnetic-only actuation in a time-varying geomagnetic field and by our previous pa\-per~\cite{KYZ2023} on a satellite with an attached flexible boom.
In this paper, we adapt the model considered in~\cite{KYZ2023} to the case of rigid-only satellite and consider an underactuated situation where the control torque is generated by magnetorquers.
For this mathematical model, we will show that the linear approximation is not controllable and prove nonlinear local controllability around a relative equilibrium under certain assumptions on the mass distribution.

\section{Mathematical model}
Consider a satellite which moves as a rigid body in a circular orbit around the Earth.
We fix a right Cartesian frame $Oxyz$ (body coordinate system -- BCS) with unit vectors $(e_1,e_2,e_3)$ assuming that $O$ is the center of mass of the satellite.
To describe the attitude of the satellite, we also consider the orbit coordinate system (OCS) as the right Cartesian frame $Ox'y'z'$ with unit vectors $(\eta_1,\eta_2,\eta_3)$ such that $\eta_1$ is orthogonal to the orbital plane, $\eta_2$ is tangent to the orbit, and $\eta_3$ points at the zenith (i.e. away from the Earth's center).

We denote by $\omega = \omega_1 e_1 + \omega_2 e_2 + \omega_3 e_3$ the absolute angular velocity vector of the satellite and
by $\omega_0 = \textrm{const}\neq 0$ the orbital rate. Then, the absolute
angular velocity of the OCS is $-\omega_0 \eta_1$, and the relative angular velocity of the satellite
with respect to the OCS is $\omega_r = \omega + \omega_0 \eta_1$.
We represent the rotation of the BCS with respect to the OCS by the attutude quaternion with the vector part $q=(q_1,q_2,q_3)$ and the scalar part $q_4$.
The equations of motion of the considered model are described by Euler's equations with respect to the angular velocity components and by the kinematic equations in quaternion form~\cite{KYZ2023,markley2014fundamentals}:
\begin{equation}\label{eq:EulerEq}
	I\dot{\omega}+\omega\times I \omega =\tau_c + 3\,\omega_0^2(\eta_3 \times I \eta_3),
\end{equation}
\begin{equation}\label{eq:QtrnEq}
\begin{aligned}
    & {\dot q}=\frac12q_4{\omega}_r + \frac12 {q} \times {\omega}_r,\\
    & \dot q_4 = - \frac12 \left<{q},{\omega}_r\right>.
    \end{aligned}
\end{equation}
We assume that $Oxyz$ are the principal axes of inertia of the satellite, so that the inertia tensor is diagonal: ${I=\textrm{diag}(I_1,I_2,I_3)}$.
The term $\tau_c$ in~\eqref{eq:EulerEq} is the torque generated by magnetorquers according to the following law:
\begin{equation}
\tau_c = \beta u \times \eta_1,\quad \beta=\textrm{const} \neq 0,
\label{magn_torque}
\end{equation}
where $u\in {\mathbb R}^3$ is the control input (magnetic moment generated by the satellite), and $\beta\eta_1$ is the Earth magnetic field vector.
Formula~\eqref{magn_torque} corresponds to the case of an equatorial orbit under the assumption of a constant geomagnetic field.
The last term in the right-hand side of~\eqref{eq:EulerEq} describes the gravity gradient torque~\cite{wertz}.

For a component-wise representation of the above equations, we express the unit vectors of the OCS by their coordinates in the BCS~\cite{MWZ2022}:
\begin{equation}\label{eq:Vec_k}
{\eta_1}=\begin{pmatrix}
q_1^2-q_2^2-q_3^2+q_4^2 \\
2(q_1 q_2-q_3 q_4)\\
2(q_1 q_3+q_2 q_4)
\end{pmatrix},\;
{\eta_2}=\begin{pmatrix}
2(q_1 q_2+q_3 q_4) \\
q_4^2-q_1^2+q_2^2-q_3^2\\
2(q_2 q_3-q_1 q_4)
\end{pmatrix},\;
{\eta_3}=
    \left(
      \begin{array}{c}
        2(q_1q_3-q_2q_4) \\
        2(q_2q_3+q_1q_4) \\
        q_3^2+q_4^2-q_1^2-q_2^2 \\
      \end{array}
    \right).
\end{equation}
Formulas~\eqref{magn_torque} and~\eqref{eq:Vec_k} allow presenting differential equations~\eqref{eq:EulerEq} and~\eqref{eq:QtrnEq} with respect to $\omega(t)$ and $q(t)$ as the following control system:
\begin{equation}\label{sys}
    \dot{x} = f(x,u),\quad x=(\omega_1, \omega_2, \omega_3, q_1, q_2, q_3)^\top\in D\subset {\mathbb R}^6,\;u=(u_1,u_2,u_3)^\top\in{\mathbb R}^3,
\end{equation}
where $D=\{(\omega_1, \omega_2, \omega_3, q_1, q_2, q_3)^\top\in{\mathbb R}^6 \,:\, q_1^2+q_2^2+q_3^2<1\}$ and
\begin{equation*}
    f(x,u) =
    \begin{pmatrix}
        \tfrac{I_2-I_3}{I_1}\omega_2\omega_3+ \tfrac{6\omega_0^2(I_2-I_3)}{I_1} (q_1q_4+q_2q_3)(2q_1^2+2q_2^2-1)+\tfrac{1}{I_1}(u_2b_3-u_3b_2)\\
        \tfrac{I_3-I_1}{I_2}\omega_1\omega_3+\tfrac{6\omega_0^2(I_3-I_1)}{I_2}(q_1q_3-q_2q_4)(2q_1^2+2q_2^2-1)+\tfrac{1}{I_2}(u_3b_1-u_1b_3)\\
        \tfrac{I_1-I_2}{I_3}\omega_1\omega_2+\tfrac{12\omega_0^2(I_1-I_2)}{I_3}(q_1 q_4 + q_2 q_3)(q_2q_4-q_1q_3)+\tfrac{1}{I_3}(u_1b_2-u_2b_1)\\
        \tfrac12 ((\omega_1+\omega_0)q_4-\omega_2q_3+\omega_3q_2)\\
        \tfrac12 ((\omega_1-\omega_0)q_3+\omega_2q_4-\omega_3q_1)\\
        \tfrac12 (-(\omega_1-\omega_0)q_2+\omega_2q_1+\omega_3q_4)
    \end{pmatrix},
\end{equation*}
\begin{equation}
q_4=\sqrt{1-q_1^2-q_2^2-q_3^2},
\label{q4eq}
\end{equation}
$$
b_1=\beta(1-2q_2^2-2q_3^2),\; b_2=2\beta(q_1 q_2-q_3 q_4),\; b_3=2\beta(q_1 q_3+q_2 q_4).
$$

We do not consider the differential equation for $\dot q_4$ in~\eqref{sys} because the attitude quaternion is normalized (${q_1^2+q_2^2+q_3^2+q_4^2=1}$), so in our case, the component $q_4$ is defined by~\eqref{q4eq} for all $x\in D$.

Thus, we have derived {\em the nonlinear control system~\eqref{sys} as a mathematical model of the considered orbiting satellite}.
This control system can be obtained from the mathematical model of~\cite{KYZ2023} by neglecting the deformation of the flexible part of the satellite.

\section{Main results}

The considered system~\eqref{sys} admits an equilibrium $(x_e,u_e)\in D\times {\mathbb R}^3$ of the form
\begin{equation*}
(x_e,u_e):\; \omega_1=-\omega_0, \omega_2=\omega_3=0,\, q_1=q_2=q_3=0,\, u_1=u_2=u_3=0,
\end{equation*}
i.e. $f(x_e,u_e)=0$. This equilibrium corresponds to the circular motion of the satellite with the constant orbital rate $\omega_0$ when the BCS aligns with the OCS.

To study the controllability property of system~\eqref{sys} near this equilibrium, we recall the notion of small-time local controllability.
For a point $\xi \in{\mathbb R}^n$ and $\varepsilon>0$, we denote the $\varepsilon$-neighborhood of $\xi$ by $B_\varepsilon(\xi)=\{\tilde
\xi\in\mathbb{R}^n:\|\tilde\xi-\xi\|<\varepsilon\}$. All the norms in this paper are Euclidean.

{\bf Definition~1.}
System~\eqref{sys} is called {\em small-time locally controllable} at the equilibrium $(x_e,u_e)$ if, for every $\varepsilon > 0$, there exists $\delta > 0$ such that, for every $x^{(0)},x^{(1)}\in B_\delta(x_e)$, there exists a function $u\in L^\infty([0,\varepsilon];\mathbb{R}^3)$ such that
\begin{equation*}
\|u(t)-u_e\|\leqslant\varepsilon\qquad\text{for all}\;t\in[0,\varepsilon],
\end{equation*}
and the solution $x(t)$ of the corresponding Cauchy problem
\begin{equation*}
    \dot x(t) = f(x(t), u(t)),\qquad x(0)=x^{(0)}
\end{equation*}
is such that
\begin{equation*}
    x(\varepsilon)=x^{(1)}.
\end{equation*}

Our main result reads as follows.

{\bf Theorem~1.}~{\em
System~\eqref{sys} is small-time locally controllable at the equilibrium $(x_e,u_e)$, provided that
$I_1\neq I_2+I_3$ and $I_2\neq I_3$.}

{\em Remark.}
From the physical viewpoint, the restrictions of Theorem~1 mean that the mass distribution is not degenerate to a flat body case ($I_1\neq I_2+I_3$), and the body is not symmetric with respect to its first principal axis of inertia ($I_2 \neq I_3$).

The proof of Theorem~1 will be given in Section~4 using Lie algebraic sufficient controllability conditions due to H.J.~Sussmann~\cite{Sus1987}.
It should be emphasized that {\em the linear approximation of system~\eqref{sys} at $(x_e,u_e)$ is not controllable}.
To show this, we first identify the components of $x$ as
$$
x_1=\omega_1,\; x_2 = \omega_2,\; x_3 = \omega_3,\; x_4 = q_1,\; x_5 = q_2,\; x_6 = q_3,
$$
and compute the Jacobian matrices $\frac{\partial f(x,u)}{\partial x}=\left(
\frac{\partial f_i(x,u)}{\partial x_j}
\right)$, $\frac{\partial f(x,u)}{\partial u}=\left(
\frac{\partial f_i(x,u)}{\partial u_j}
\right)$ at $(x_e,u_e)$:
\begin{equation*}
A=\left.\frac{\partial f(x,u)}{\partial x}\right|_{(x,u)=(x_e,u_e)}=
\begin{pmatrix}
0 & 0 & 0 & \tfrac{6\omega_0^2(I_3-I_2)}{I_1} & 0 & 0
\\ 0 & 0 & \tfrac{I_1-I_3}{I_2}\omega_0 & 0 & \tfrac{6\omega_0^2(I_3-I_1)}{I_2} & 0
\\ 0 & \tfrac{I_2-I_1}{I_3}\omega_0 & 0 & 0 & 0 & 0
\\ \tfrac12 & 0 & 0 & 0 & 0 & 0
\\ 0 & \tfrac12 & 0 & 0 & 0 & -\omega_0
\\ 0 & 0 & \tfrac12 & 0 & \omega_0 & 0
\end{pmatrix},
\end{equation*}
\begin{equation*}
B=\left.\frac{\partial f(x,u)}{\partial u}\right|_{(x,u)=(x_e,u_e)}=
\begin{pmatrix}
   0 & 0 & 0
\\ 0 & 0 & \tfrac{\beta}{I_2}
\\ 0 & -\tfrac{\beta}{I_3} & 0
\\ 0 & 0 & 0
\\ 0 & 0 & 0
\\ 0 & 0 & 0
\end{pmatrix}.
\end{equation*}
Then the linearization of~\eqref{sys} at the considered equilibrium reads as
\begin{equation}\label{linear}
    \dot \xi = A \xi +B v,\quad \xi \in {\mathbb R}^6,\; v\in {\mathbb R}^3.
\end{equation}

The rank of the block matrix $[B,AB,A^2B,A^3B,A^4B,A^5B]$ is at most $4$; therefore, Kalman's controllability condition for system~\eqref{linear} fails. This means that the linear test cannot be applied to the local controllability analysis of system~\eqref{sys} near the equilibrium $(x_e,u_e)$.

\section{Small-time local controllability: Proof of Theorem~1}


For further analysis, we rewrite system~\eqref{sys} in the following control-affine form:

\begin{equation}\label{eq3.1}
    \dot{x} = f_0({x})+\sum_{i=1}^3 u_i f_i({x}),
\end{equation}
where the smooth vector fields $f_i:D\to {\mathbb R}^6$ are
\begin{equation*}
    f_0({x}) =
    \begin{pmatrix}
        \tfrac{I_2-I_3}{I_1}\omega_2\omega_3+ \tfrac{6\omega_0^2(I_2-I_3)}{I_1} (q_1q_4+q_2q_3)(2q_1^2+2q_2^2-1)\\
        \tfrac{I_3-I_1}{I_2}\omega_1\omega_3+\tfrac{6\omega_0^2(I_3-I_1)}{I_2}(q_1q_3-q_2q_4)(2q_1^2+2q_2^2-1)\\
        \tfrac{I_1-I_2}{I_3}\omega_1\omega_2+\tfrac{12\omega_0^2(I_1-I_2)}{I_3}(q_1 q_4 + q_2 q_3)(q_2q_4-q_1q_3)\\
        \tfrac12 ((\omega_1+\omega_0)q_4-\omega_2q_3+\omega_3q_2)\\
        \tfrac12 ((\omega_1-\omega_0)q_3+\omega_2q_4-\omega_3q_1)\\
        \tfrac12 (-(\omega_1-\omega_0)q_2+\omega_2q_1+\omega_3q_4)
    \end{pmatrix},
\end{equation*}

\begin{equation*}
     f_1(x)=
    \begin{pmatrix}
        0\\
        -\tfrac{2\beta(q_1 q_3+q_2 q_4)}{I_2}\\
        \tfrac{2\beta(q_1 q_2-q_3 q_4)}{I_3}\\
        0\\
        0\\
        0
    \end{pmatrix}, \,
    f_2(x) =
    \begin{pmatrix}
        \tfrac{2\beta(q_1 q_3+q_2 q_4)}{I_1}\\
        0\\
        \tfrac{\beta(2q_2^2+2q_3^2-1)}{I_3}\\
        0\\
        0\\
        0
    \end{pmatrix},\,
    f_3(x) =
    \begin{pmatrix}
        \tfrac{2\beta(q_3 q_4-q_1 q_2)}{I_1}\\
        \tfrac{\beta(1-2q_2^2-2q_3^2)}{I_2}\\
        0\\
        0\\
        0\\
        0
    \end{pmatrix},
\end{equation*}
and the component $q_4$ is defined by~\eqref{q4eq}.

We denote the above family of vector fields by $\mathcal{F}=\{f_0,f_1,...,f_m\}$, $m=3$. For $f_i,f_j\in \mathcal{F}$, the Lie bracket of $f_i$ and $f_j$ is
$$
[f_i,f_j](x) = \frac{\partial f_j(x)}{\partial x} f_i(x) - \frac{\partial f_i(x)}{\partial x} f_j(x),
$$
where $\frac{\partial f_i(x)}{\partial x}$ is the Jacobian matrix of $f_i(x)$.
The {\em Lie algebra} of vector fields generated by $\mathcal{F}$ will be denoted by  $\mathcal{L}(\mathcal{F})$ (for more clarifications, see \cite[Section 3.2]{Cor2007}).
To prove Theorem~1, we will use Sussmann's controllability condition~\cite{Sus1987}.

{\bf Definition~2.}
The control affine-system~\eqref{eq3.1} satisfies the \textit{Lie algebra rank condition} at $(x_e, u_e)$ if
\begin{equation}\label{rank_cond}
\{g(x_e):g\in\mathcal{L}(\mathcal{F})\}=\mathbb{R}^6.
\end{equation}

Now we denote by $Br(\mathcal{F})$ the set of iterated Lie brackets of $\mathcal{F}$ and, for an $h\in Br(\mathcal{F})$, we denote by $\delta_i(h)$ the number of times that $f_i$ appears in $h$.

{\bf Definition~3.}
A control-affine system $\dot x = f_0(x)+\sum_{i=1}^m u_i f_i(x)$ satisfies the \textit{Suss\-mann condition} $S(\theta)$ for some $\theta\in[0,+\infty)$ if, for any $h\in Br(\mathcal{F})$ with $\delta_0(h)$ odd and $\delta_i(h)$ even for all $i=\overline{1,m}$, the following condition holds:
\begin{equation}\label{Sus}
   h(x_e) \in D_\theta(h):=\text{span}\{g(x_e):g\in Br(\mathcal{F}) \text{ and } \delta(g)<\delta(h)\},
\end{equation}
where
\begin{equation*}
   \delta(g):=\theta\delta_0(g)+\sum_{i=1}^m \delta_i(g).
\end{equation*}

We refer to the following sufficient conditions for small-time local controllability due to H.J.~Sussmann~\cite{Sus1987}.

{\bf Theorem~2.}~{\em
If system~\eqref{eq3.1} satisfies the Lie algebra rank condition at an equilibrium $(x_e, u_e)$ and the Sussmann condition $S(\theta)$ with some $\theta\in[0,1]$, then system~\eqref{eq3.1} is small-time locally controllable.}

\begin{proof}[Proof of Theorem~1]
We check the conditions of Theorem~2 by straightforward computation, in a similar way as it was done in \cite[Theorem~2]{Ker1995}.
It is common to call $h\in Br(f)$ a ``bad'' bracket if it satisfies $\delta_0(h)$ odd and $\delta_i(h)$ is even for every $i=\overline{1,3}$, and a ``good'' bracket if this verification fails (cf.~\cite{Cor2007,Ag2024}). The main task is to find six ``good'' brackets $g_i$ which span $\mathbb{R}^6$ at $x=x_e$, and then analyze all possible cases for $h$ in terms of condition~\eqref{Sus}.

As a set of ``good'' brackets, we take
\begin{equation}
\begin{aligned}
&g_1=[f_0,[f_3,[f_0,[f_0,f_2]]]],\;g_2=f_3,\;g_3=f_2,\\
&g_4=[f_2,[f_0,[f_0,f_3]]],\;g_5=[f_0,f_3],\;g_6=[f_0,f_2].
\end{aligned}
\label{good_brackets}
\end{equation}
The evaluation of these vector fields at the equilibrium yields
\begin{equation*}
\begin{aligned}
& g_1(x_e)=\tfrac{3\omega_0^2\beta^2(I_1-I_2-I_3)(I_2-I_3)}{I_1^2I_2I_3}\,(1,0,0,0,0,0)^\top,
\\& g_2(x_e)=\tfrac{\beta}{I_2}\,(0,1,0,0,0,0)^\top,
\\& g_3(x_e)=-\tfrac{\beta}{I_3}\,(0,0,1,0,0,0)^\top,
\\& g_4(x_e)=\tfrac{\beta^2(I_2+I_3-I_1)}{2I_1I_2I_3}\,(0,0,0,1,0,0)^\top,
\\& g_5(x_e)=\left(0,0,\tfrac{\omega_0\beta(I_1-I_2)}{I_2I_3},0,-\tfrac{\beta}{2I_2},0\right)^\top,
\\& g_6(x_e)=\left(0,\tfrac{\omega_0\beta(I_1-I_3)}{I_2I_3},0,0,0,\tfrac{\beta}{2I_3}\right)^\top.
\end{aligned}
\end{equation*}
Under the conditions of Theorem~1, it is easy to see that
$$
\text{span}\{g_1(x_e),g_2(x_e),g_3(x_e),g_4(x_e),g_5(x_e),g_6(x_e)\}=\mathbb{R}^6,
$$ which proves that the Lie algebra rank condition~\eqref{rank_cond} holds.

Taking $\theta=\tfrac{1}{2}$, we have
\begin{equation*}
\delta(g_i)\leqslant \frac{7}{2}\quad\text{for all}\; i=\overline{1,6}.
\end{equation*}
Now for ``bad'' brackets $h\in Br(f)$ with $\delta(h)>\frac{7}{2}$, Sussmann's condition is obviously satisfied. If $\delta(h)\leqslant\frac{7}{2}$, then we have three possibilities:
\begin{enumerate}
    \item[1)] $\delta_0(h)\leqslant7$, $\delta_i(h)=0$ for all $ i=\overline{1,3}$. In this case, the corresponding iterated Lie brackets vanish at $x_e$.

    \item[2)]  In the case where $\delta_0(h)=1$ and $\delta_i(h)=2$ for one of $i\in\{1,2,3\}$, we have
    \begin{equation*}
        [f_i,[f_i,f_0]] \equiv [f_0,[f_i,f_i]] \equiv 0\quad \text{for all}\; i=\overline{1,3}.
    \end{equation*}

    \item[3)]  In the case where $\delta_0(h)=3$ and $\delta_i(h)=2$ for one of $i\in\{1,2,3\}$, the only nonzero brackets at $x_e$ are:
    \begin{equation*}
    \begin{aligned}
        &h_1=[f_2,[f_0,[f_0,[f_0,f_2]]]]=-\tfrac{\omega_0 \beta^2(I_1^2 - (I_2 + 2 I_3)I_1 + 2I_2 I_3 + I_3^2)}{I_1I_2I_3^2}\,(0,0,0,1,0,0)^\top,
        \\& h_2=[f_3,[f_0,[f_0,[f_0,f_3]]]]=-\tfrac{\omega_0 \beta^2(I_1^2 -(2I_2 + I_3)I_1 + I_2 (I_2 + 2 I_3))}{I_1I_2^2I_3}\,(0,0,0,1,0,0)^\top.
    \end{aligned}
    \end{equation*}
As $\delta(g_4)=3<\delta(h_i)$ for $i=1,2$, we state that $h_i\in D_{\frac{1}{2}}(h_i)$ for $i=1,2$.
\end{enumerate}

Hence, Sussmann's condition holds which concludes the proof of Theorem~1.
\end{proof}

\section{Conclusion}

A particular outcome of the proof of Theorem~1 provides the set of ``good'' Lie brackets~\eqref{good_brackets} for control system~\eqref{eq3.1}.
According to the available stabilizability results for essentially nonlinear control systems~\cite{ZG,GZ},
the structure of iterated brackets in the controllability condition is crucial for defining the stabilizing controllers.
Control design issues for system~\eqref{eq3.1} under the bracket-generating condition, with the brackets of the form~\eqref{good_brackets}, are left as prospective directions for future research.

\bibliographystyle{unsrt}

\end{document}